\documentclass[10pt,a4paper,final]{article}
\usepackage[latin1]{inputenc}
\usepackage{amsmath}
\usepackage{amsfonts}
\usepackage{amssymb}
\usepackage{makeidx}
\usepackage{graphicx}
\usepackage{hyperref}
\usepackage[notcite,notref,color]{showkeys}

\usepackage{todo}
\usepackage{amssymb}
\usepackage{amsmath}
\usepackage{latexsym}
\usepackage{amssymb}
\usepackage{amsthm}
\usepackage{cite}
\usepackage{mathrsfs}
\usepackage{bbm}
\usepackage{bbold}
\usepackage{subfigure}
\usepackage{float}
\usepackage{epsfig}
\usepackage[numbers]{natbib}
\usepackage{mathtools}
\usepackage{listings}
\usepackage{pgf}
\usepackage{listings}
\usepackage{epstopdf}

\theoremstyle{plain} 
\newtheorem{thm}{Theorem}[section] 
\newtheorem{prop}[thm]{Proposition}

\newtheorem{lem}[thm]{Lemma}

\theoremstyle{definition} 
\newtheorem{defn}{Definition}[section]

\theoremstyle{remark} 
\newtheorem{oss}{Remark}

\title{An elementary approach to component sizes in critical random graphs}
\author{Umberto De Ambroggio\thanks{University of Bath, Department of Mathematical Sciences - umbidea@gmail.com}}
\begin{document}
	\maketitle
\begin{abstract}
	In this article we introduce a simple tool to derive polynomial upper bounds for the probability of observing unusually large maximal components in some models of random graphs when considered at criticality. Specifically, we apply our method to a model of random intersection graph, a random graph obtained through $p$-bond percolation on a general $d$-regular graph, and a model of inhomogeneous random graph.
	\vskip0.1cm
	MSC class: 60C05, 05C80.
\end{abstract}

\section{Introduction} 
The purpose of this paper is to introduce an elementary tool to obtain simple upper bounds for the probability of observing unusually large maximal components in some important models of random graphs at criticality.

Given any (undirected) graph $\mathbb{G}=(V,E)$ and vertices $v,u\in V$, we write $v\sim u$ if the edge $\{v,u\}$ is present in $\mathbb{G}$ and say that vertices $u$ and $v$ are neighbours. We write $v\leftrightarrow u$ if there exists a path of occupied edges connecting vertices $v$ and $u$ and we adopt the convention that $v\leftrightarrow v$ for every $v\in V$. We denote by $\mathcal{C}(v)= \{u\in V:v\leftrightarrow u\}$ the component (or cluster) of vertex $v\in V$ and its size by $|\mathcal{C}(v)|$. Moreover, we define a largest component $\mathcal{C}_{\max}$ to be any cluster $\mathcal{C}(v)$ for which $|\mathcal{C}(v)|$ is maximal, so that $|\mathcal{C}_{\max}|=\max_{v\in V}|\mathcal{C}(v)|$.

The Erd\H{o}s-R\'enyi random graph on $[n]= \{1,\dots,n\}$, denoted by $\mathbb{G}(n,p)$, is the random graph obtained from the complete graph on $n$ vertices by independently retaining each edge with probability $p\in [0,1]$ and deleting it with probability $1-p$. 

One of the most surprising aspects of this model is that when $p$ is of the form $p=p(n)=\gamma/n$, then the $\mathbb{G}(n,p)$ random graph undergoes a \textit{phase transition} as $\gamma$ passes 1. Specifically, if $\gamma < 1$, then $|\mathcal{C}_{\max}|$ is of order $\log(n)$; if $\gamma=1$, then $|\mathcal{C}_{\max}|$ is of order $n^{2/3}$; and if $\gamma > 1$, then $|\mathcal{C}_{\max}|$ is of order $n$. See for instance the books \cite{bollobas_book},\cite{janson_et_al:random_graphs} or \cite{remco:random_graphs} for proofs of these statements and other interesting properties of this model. See also Krivelevich and Sudakov \cite{KrivSud} for a simple proof of the phase transition in $\mathbb{G}(n,p)$.

In \cite{de_ambroggio_roberts1} De Ambroggio and Roberts introduced a ballot-type result (Lemma \ref{ballot} below) to provide a new, purely probabilistic proof of the fact that in the $\mathbb{G}(n,p)$ model considered in the so-called \textit{critical window}, i.e. when $p$ is of the form $p=p(n)=n^{-1}+\lambda n^{-4/3}$, the probability of observing a maximal cluster of size larger than $An^{2/3}$ tends to zero (as $n\rightarrow \infty$) exponentially fast in $A$. More precisely, they proved that, for large enough $n$, 
\begin{equation}\label{erd}
\frac{c}{A^{3/2}}e^{-\frac{A^3}{8}+\frac{\lambda A^2}{2}-\frac{\lambda^2A}{2}}\leq \mathbb{P}(|\mathcal{C}_{\max}|>An^{2/3})\leq \frac{c'}{A^{3/2}}e^{-\frac{A^3}{8}+\frac{\lambda A^2}{2}-\frac{\lambda^2A}{2}},
\end{equation} 
where $0<c\leq c'$ are two finite constants, thus showing that in the near-critical $\mathbb{G}(n,p)$ model the number of vertices contained in the maximal component is unlikely to be much larger than $n^{2/3}$.

We remark that the correct asymptotic for $\mathbb{P}(|\mathcal{C}_{\max}|>An^{2/3})$ in this critical model was obtained first by Pittel \cite{pittel:largest_cpt_rg} (whose paper is partially based on an earlier article by Luczak, Pittel and Wierman \cite{luczak_et_al:structure_RG}) and more recently by Roberts \cite{roberts:component_ER}. We also mention that Nachmias and Peres \cite{nachmias_peres:CRG_mgs} used a general martingale argument to establish an exponential upper bound for the probability in (\ref{erd}), but their bound is not optimal.

The purpose of this work is to show that part of the argument used in \cite{de_ambroggio_roberts1} to prove the upper bound in (\ref{erd}) is quite general and can be used to obtain, in a surprisingly simple way, polynomial upper bounds for $\mathbb{P}(|\mathcal{C}_{\max}|>k)$ in different models of random graphs when considered at criticality. Specifically, we apply our method to three different models, namely a model of random intersection graph, a random graph obtained through $p$-bond percolation on a general $d$-regular graph, and a model of inhomogeneous random graph, and we show that $|\mathcal{C}_{\max}|$ is unlikely to be much larger than $n^{2/3}$ in these models. In this sense, these random graphs exhibit a similar critical behaviour.

\section{Results}
In order to better understand the statement of our main result (Theorem \ref{Umbi} below), we first need to recall the definition of an \textit{exploration process}, which is an algorithmic procedure used to reveal the components of a given graph; see e.g. \cite{de_ambroggio_roberts1}, \cite{nachmias_peres:CRG_mgs}, \cite{nachmiasperes_criticalperc}, \cite{roberts:component_ER} and references therein. As we will see in a moment, when the graph under investigation is \textit{random} such exploration process reduces the study of component sizes to the analysis of the trajectory of a random process, which looks like (but it is not quite) a random walk.

Let $\mathbb{G}=([n],E)$ be any (undirected) graph, and let $V_n$ be a vertex selected uniformly at random from $[n]$. During the exploration of $\mathcal{C}(V_n)$, each vertex will be either \textit{active}, \textit{explored} or \textit{unseen}, and its status will change during the course of the exploration. At each step $t\in \{0\}\cup [n]$ of the algorithm, the number of explored vertices will be $t$ whereas the number of active vertices will be denoted by $Y_t$. At time $t=0$, if $V_n$ is an isolated vertex we stop the procedure; otherwise, there exists some vertex $u\in [n]\setminus {V_n}$ with $\{V_n,u\}\in E$. In this case, vertices $V_n$ and $u$ are declared active, whereas all other vertices are declared unseen (so that $Y_0=2$). At each step $t\in [n]$ of the algorithm we proceed as follows.
\begin{itemize}
	\item [(a)] If $Y_{t-1}>0$, we let $u_t$ be the active vertex with the smallest label. We reveal all unseen neighbours of $u_t$ in $\mathbb{G}$ and change the status of these vertices to active. Then, we set $u_t$ itself explored. 
	\item [(b)] If $Y_{t-1}=0$, we let $u_t$ be the unseen vertex with the smallest label, and:
	\begin{itemize}
		\item [(b.1)] if $u_t$ is isolated, we halt the procedure;
		\item  [(b.2)] otherwise, there is at least one unseen vertex $v$ such that $\{u_t,v\}\in E$, and we declare active both $u_t$ and $v$; then we continue with step (a).
	\end{itemize}
\end{itemize}
Denoting by $\eta_t$ the number of unseen vertices in $\mathbb{G}$ which become active at step $t$ of the exploration process, we see that:
\begin{itemize}
	\item [(i)] if $Y_{t-1}>0$ then $Y_t=Y_{t-1}+\eta_t-1$;
	\item [(ii)] if $Y_{t-1}=0$ then $Y_t=\eta_t$.
\end{itemize}
\begin{oss}
	We remark that our description of an exploration process is slightly different from those provided e.g. in \cite{de_ambroggio_roberts1}, \cite{nachmias_peres:CRG_mgs} and \cite{roberts:component_ER}. Indeed, in our setting the algorithm is actually run only in the case where $V_n$ is not an isolated vertex, and the exploration starts from two active vertices and not from one active vertex, as it usually happens. Moreover, whenever $Y_{t-1}=0$ and $u_t$ is \textit{not} isolated, we first reveal one of the neighbours of $u_t$ before proceeding with the exploration. These small modifications will be particularly useful in one of our applications.
\end{oss}
Now let $\mathbb{G}=([n],E)$ be any (undirected) random graph, and use the above algorithm to reveal the components of $\mathbb{G}$. It is clear that now the $\eta_i$ are random variables. Observe that, given any $k\in \mathbb{N}=\{1,2,\dots\}$, if $V_n$ is an isolated vertex then $|\mathcal{C}(V_n)|=1$ and hence, in particular, we can't have $|\mathcal{C}(V_n)|>k$ (as $k\geq 1$). On the other hand, if $V_n$ is not isolated then $|\mathcal{C}(V_n)|>k$ implies that $Y_t=2+\sum_{i=1}^{t}(\eta_{i}-1)>0$ for all $t\in [k]$. Therefore we can write
\begin{align}\label{mainid}
\mathbb{P}(|\mathcal{C}(V_n)|>k)\leq \mathbb{P}\left(2+\sum_{i=1}^{t}(\eta_{i}-1)>0 \hspace{0.15cm}\forall t\in [k]\right).
\end{align}
Note that the $\eta_i$ are not independent and, moreover, they have different distributions (one of the reasons is that the number of unseen vertices in the graph decreases during the course of the exploration). Therefore $Y_t$ does not define a random walk. 

In order to bound from above the probability in (\ref{mainid}), the idea is to produce a sequence of independent and identically distributed (i.i.d.)~random variables $X_i$, \textit{bigger} than the $\eta_{i}$, that allow us to replace the probability on the right hand side of (\ref{mainid}) with the probability that a random walk (started at 2) stays positive up to time $k$. 

In some random graphs this is an immediate consequence of the model construction, while in other instances one needs more care in order to produce these $X_i$.

Here is our main result.
\begin{thm}\label{Umbi}
	Let $\mathbb{G}=([n],E)$ be any (undirected) random graph. Suppose that there exists a sequence of i.i.d.~random variables $(X_i)_{i\geq 1}$ taking values in $\mathbb{N}_0$, such that the distribution of $X_1$ may depend on $n$, satisfying: 
	\begin{itemize}
		\item [(i)] $\mathbb{P}(X_1=3)\geq c$ for all sufficiently large $n$, for some constant $c>0$;
		\item [(ii)] $\mathbb{P}(|\mathcal{C}(V_n)|>k)\leq \mathbb{P}\left(2+\sum_{i=1}^{t}(X_i-1)>0\hspace{0.15cm}\forall t\in [k]\right)$ for every $k=k(n)\in \mathbb{N}$;
		\item [(iii)] there exist $\delta,\rho=\rho(n)>0$ and $\epsilon=\epsilon(n)\geq 0$ with $\epsilon_n\rightarrow 0$ (as $n\rightarrow \infty$) such that $\mathbb{E}\left(e^{rX_1}\right)\leq e^{r(1+\epsilon)+r^2\delta}$ for every $r\in (0,\rho)$ and all sufficiently large $n$.
	\end{itemize}
	Suppose that $k=k(n)\in \mathbb{N}$ satisfies $\epsilon \sqrt{k}\leq c$ and $\rho \sqrt{k}\geq 1$ for all large enough $n$, for some finite constant $c>0$. Then
	\begin{equation}\label{bound}
	\mathbb{P}(|\mathcal{C}_{\max}|>k)\leq \frac{C}{\mathbb{P}(X_1=3)}\frac{n}{k^{3/2}}
	\end{equation}
	for all sufficiently large $n$, where $C>0$ is a finite positive constant which depends solely on $\delta$ and $c$.
\end{thm}

\begin{oss}\label{rem1}
	In all our applications the probability $\mathbb{P}(X_1=3)$ is bounded away from zero, so that if $k=k(n)= \lceil An^{2/3} \rceil$ satisfies the two assumptions in the statement of the theorem then $\mathbb{P}(|\mathcal{C}_{\max}|> \lceil An^{2/3} \rceil)$ is indeed $O(A^{-3/2})$.
\end{oss}
\begin{oss}
	We note that condition $(iii)$ in Theorem \ref{Umbi} might be stated in different (possibly more general) terms, but we decided to state it in this way because of its simplicity to be verified, as shown in our applications.
\end{oss}
Our claim that the approach introduced in \cite{de_ambroggio_roberts1} is robust and that Proposition \ref{Umbi} leads to simple upper bounds for $\mathbb{P}(|\mathcal{C}_{\max}|>k)$ in several models of random graphs at criticality is justified in Sub-sections 2.1, 2.2 and 2.3 below, where we apply Proposition \ref{Umbi} to obtain polynomial upper bounds for the above probability in three particular models of random graphs.

We remark that our methodology does not lead to upper bounds for the probability of observing unusually \textit{small} largest components; in this direction the martingale argument introduced by Nachmias and Peres in \cite{nachmias_peres:CRG_mgs} seems to be more robust and adaptable to different models of random graphs.


\subsection{Critical random intersection graph.}
Our first application of Theorem \ref{Umbi} involves a model of random intersection graph; for an introduction to this class of models, we refer the reader to \cite{frieze_karo?ski_2015}.

Here we are interested in the random graph described by Lageras and Lindholm \cite{laglind2008}. Such a random graph, denoted by $\mathbb{G}(n,m,p)$, with a set of vertices $V=\{v_i:i\in [n]\}$ and a set of edges $E$, is constructed from a bipartite graph $\mathbb{B}(n,m,p)$ with two sets of vertices: $A=\{a_j:j\in [m]\}$, which we call the set of \textit{auxiliary vertices}, and $V$ (that is, the vertex set of $\mathbb{G}(n,m,p)$). Edges in $\mathbb{B}(n,m,p)$ between vertices and auxiliary vertices are present independently with probability $p\in [0,1]$. Two distinct vertices $v_i$ and $v_j$ are neighbours in $\mathbb{G}(n,m,p)$ (i.e. $\{v_i,v_j\}\in E$) if and only if there exists at least one $a_k\in A$ such that both edges $\{a_k,v_i\}$ and $\{a_k,v_j\}$ are present in the bipartite graph $\mathbb{B}(n,m,p)$. 

We are interested in the case where $p=p(n)=\gamma/n^{(1+\alpha)/2}$ and $m=m(n)=
\lfloor \beta n \rfloor$, where $\alpha,\beta,\gamma> 0$ are fixed parameters of the model.   

Stark \cite{starkRIG} has shown that the vertex degree distribution (i.e. the distribution of the degree of a vertex selected uniformly at random) is highly dependent on the value of $\alpha$. However, as shown by Deijfen and Kets \cite{deijfen_kets_2009}, the clustering is controllable only when $\alpha=1$. 

The component structure of the graph is studied for $\alpha\neq 1$,$\gamma>0$ and $\beta=1$ by Bherisch \cite{beh2007}, whereas it is studied for $\alpha=1$ and $\beta, \gamma>0$ in \cite{laglind2008}. Specifically, Lageras and Lindholm \cite{laglind2008} proved that the $\mathbb{G}(n,m,p)$ model undergoes a phase transition as $\beta \gamma^2$ passes $1$. Indeed, setting $\mu= \beta \gamma^2$, they proved that if $\mu<1$ (sub-critical case) then with probability tending to one there is no component in $\mathbb{G}(n,m,p)$ with more than $O(\log(n))$ vertices, while if $\mu >1$ (super-critical case) then, with probability tending to one, there exists a unique giant component of size $n\delta$ where $\delta\in (0,1)$, and the size of the second largest component is at most of order $\log(n)$.

By means of Theorem \ref{Umbi} we show that, in the \textit{critical} case $\mu=1$, it is unlikely for the largest component to contain more than $n^{2/3}$ vertices.

\begin{prop}\label{CRIG}
	Let $\mathbb{G}(n,m,p)$ be the random intersection graph described above. Let $m= \lfloor \beta n \rfloor$, $p= \gamma/n$ and $\mu= \beta \gamma^2$. If $\mu=1$ then, given any $A>1$, when $n$ is sufficiently large we have that
	\begin{equation}\label{smallres}
	\mathbb{P}(|\mathcal{C}_{\max}|> \lceil A n^{2/3} \rceil)\leq \frac{c_1}{A^{3/2}},
	\end{equation}
	where $c_1$ is a finite constant which depends solely on $\gamma$ and $\beta$.
\end{prop}

\subsection{Critical $p$-bond percolation on $d$-regular graph.}
In this section we consider a second application of Theorem \ref{Umbi}. Here we analyse a random graph $\mathbb{G}_p$ obtained through $p$-bond percolation on a general $d$-regular graph.

In \cite{nachmiasperes_criticalperc} Nachmias and Peres adapted the martingale method developed by the same authors in \cite{nachmias_peres:CRG_mgs} to prove that, for any $d\geq 3$, when $p\leq (d-1)^{-1}$ then
\begin{equation}\label{peres}
\mathbb{P}(|\mathcal{C}_{\max}|> \lceil An^{2/3}\rceil) \leq \frac{8}{A^{3/2}};
\end{equation}
see Proposition $1.2$ in \cite{nachmiasperes_criticalperc}. For a random regular graph $\mathbb{G}(n,d,p)$ they were also able to sharpen the upper bound in (\ref{peres}) and to prove a corresponding lower bound. (The $\mathbb{G}(n,d,p)$ random graph is obtained by the following two-step procedure: first we draw uniformly at random a graph from the set of all simple $d$-regular graphs on $[n]$, and then we retain each edge independently with probability $p$ and delete it with probability $1-p$.) Specifically, in Theorem 2 of \cite{nachmiasperes_criticalperc} it is shown that, when $p$ is of the form
\begin{equation}\label{criticalp}
p=p(n,d)=(1+\lambda n^{-1/3})/(d-1) \hspace{0.2cm}(\lambda \in \mathbb{R})
\end{equation}
and $d\geq 3$ is \textit{fixed}, then there are constants $C_1,C_2\in (0,\infty)$ depending on $\lambda$ and $d$ such that, for every $A>0$ and all $n$, $\mathbb{P}(|\mathcal{C}_{\max}|>An^{2/3})\leq A^{-1}C_1 e^{-C_2 A^3}$. In \cite{nachmiasperes_criticalperc} it is also shown that there exists a constant $C_3\in (0,\infty)$ (also depending on $\lambda$ and $d$) such that, for large enough $A$ and all $n$, then $\mathbb{P}(|\mathcal{C}_{\max}|<\lceil A^{-1} n^{2/3}\rceil)\leq C_3A^{-1/2}$, thus proving that the size of $|\mathcal{C}_{\max}|$ is indeed of order $n^{2/3}$ in this model when considered at criticality.

We remark that in \cite{nachmiasperes_criticalperc}, the parameter $d$ is not allowed to depend on $n$. The problem of determining the size of $|\mathcal{C}_{\max}|$ in the critical $\mathbb{G}(n,d,p)$ model when $d=d(n)$ depends on $n$ has been investigated by Joos and Perarnau \cite{joosper}, where the authors proved (among many other things) that for any $d\in \{3,\dots,n-1\}$ and when $p$ is of the form (\ref{criticalp}), then for all sufficiently large $n$ and $A=A(\lambda)$ we have that $\mathbb{P}(|\mathcal{C}_{\max}|\notin [A^{-1}n^{2/3},An^{2/3}])\leq 20/\sqrt{A}$.

Our goal here is to show that, by means of Theorem \ref{Umbi}, we can recover (up to a multiplicative constant) the bound in (\ref{peres}), in a very simple way.
\begin{prop}\label{CBPRRG}
	Let $\mathbb{G}$ be a $d$-regular graph, $d> 3$, and denote by $\mathbb{G}_p$ the random graph obtained by bond percolation on $\mathbb{G}$ with probability $p$. If $p\leq 1/(d-1)$ then, given any $A>1$, when $n$ is sufficiently large we have that
	\begin{equation}\label{ff}
	\mathbb{P}(|\mathcal{C}_{\max}|> \lceil An^{2/3}\rceil)\leq \frac{c_2}{A^{3/2}},
	\end{equation}
	for some finite positive constant $c_2$ which depends solely on $d$.
\end{prop}
\begin{oss}
	The requirement $d>3$ is needed because the i.i.d.~random variables $(X_i)_{i\geq 1}$ that dominate the $\eta_i$ in the exploration process satisfy
	\begin{equation}
	\mathbb{P}(X_1=3)=\frac{(d-2)(d-3)}{6(d-1)^2}\left(1-\frac{1}{d-1}\right)^{d-4};
	\end{equation}
	see Subsection 3.3. Hence, if we want $\mathbb{P}(X_1=3)>0$ (a condition required in Theorem \ref{Umbi}), we do need $d>3$. 
\end{oss}

\subsection{Critical inhomogeneous random graph.}
In this section we discuss our final application of Theorem \ref{Umbi}. In the random graph model that we investigate here, the $n$ vertices are endowed with weights, and edges between pair of vertices are placed independently with probabilities moderated by such weights.

Specifically, let $\mathbf{w}=(w_i)_{i\in [n]}$ be a sequence of positive real numbers, which we call the sequence of \textit{vertex weights}; we think of $w_i$ as the weight assigned to vertex $i\in [n]$. Define $l_n= \sum_{i\in[n]}w_i$, the sum of all weights.

We consider the so-called Norros-Reittu random graph \cite{norros_reittu_2006} as described by Van der Hofstad \cite{hofstad_critical_behaviour}. This is an inhomogeneous random graph, that we denote by $NR_n(\mathbf{w})$, in which the probability that the edge $\{i,j\}$ is present in $NR_n(\mathbf{w})$ (for $1\leq i<j\leq n$) is given by
\begin{equation}\label{pij}
p_{ij}^{NR}= \mathbb{P}(\{i,j\}\in E(NR_n(\mathbf{w})))=  1-e^{-w_jw_j/l_n}, \hspace{0.15cm}
\end{equation}
and edges are present independently. 

Inhomogeneous random graphs were studied extensively by Bollobas, Janson and Riordan \cite{boll_janson_riord_inhom}. As explained by Janson \cite{janson_contiguity} and further remarked by Van der Hofstad \cite{hofstad_critical_behaviour}, the $NR_n(\mathbf{w})$ random graph is closely related to the models studied by Chung and Lu \cite{chung2002connected, chung2004average, chung2006volume} and Norros and Reittu \cite{norros_reittu_2006}, so that the results proved for the $NR_n(\mathbf{w})$ random graph apply as well to these other models.

Other models of inhomogeneous random graphs have been studied more recently by Penrose \cite{penrose_2018} and by Kang, Pachon and Rodriguez \cite{kang2018evolution}.

It is clear that the topology of the $NR_n(\mathbf{w})$ model depends on the choice of the sequence $\mathbf{w}$, which we now specify.

Let $F:(0,\infty)\mapsto [0,1]$ be a distribution function, and define
\begin{equation*}\label{F}
\hspace{0.3cm}[1-F]^{-1}(u)= \inf\{s:[1-F(s)]\leq u\}, \quad u\in(0,1).
\end{equation*}
By convention, we set $[1-F]^{-1}(1)= 0$. We construct the weights as in \cite{hofstad_critical_behaviour}, namely we set
\begin{equation}\label{wj}
w_j= [1-F]^{-1}(j/n),\hspace{0.2cm}j\in [n].
\end{equation}
In \cite{boll_janson_riord_inhom} (Theorem 3.13) it has been shown that in the $NR_n(\mathbf{w})$ random graph with vertex weights as in (\ref{wj}), the proportion of vertices having degree $k\geq 0$, denoted by $N_k$, converges in probability (as $n \rightarrow \infty$) to 
\[p_k= \mathbb{E}\left(e^{-W}\frac{W^k}{k!}\right),\] 
where $W$ is a random variable taking values in $(0,\infty)$ with distribution function $F$. The limiting sequence $(p_k)_{k\geq 0}$ is a so-called \textit{mixed Poisson distribution with mixing distribution $F$}. 

In order to describe the phase transition for the size of the largest component, define
\begin{equation}\label{nu}
\nu=\mathbb{E}(W^2)/\mathbb{E}(W).
\end{equation}
As explained by Van der Hofstad \cite{hofstad_critical_behaviour} (see also \cite{de_ambroggio_pachon:upper_bounds_inhom_RGs}), this (positive) real number corresponds to the asymptotic mean of the offspring distribution in a branching process approximation of the exploration of $\mathcal{C}(V_n)$.

In \cite{boll_janson_riord_inhom}(Theorem 3.1) it is shown that the graph undergoes a phase transition as $\nu$ passes $1$. In particular, if $\nu>1$, the largest component contains a positive proportion of the total number of vertices, whereas if $\nu \leq 1$ the largest component contains a vanishing proportion of vertices.

Van der Hofstad \cite{hofstad_critical_behaviour} provided a complete picture of the component structure in the critical $NR_n(\textbf{w})$ model. Specifically, he proved that in the case where  
\begin{equation}\label{five}
\lim_{x\rightarrow \infty}x^{\tau-1}(1-F(x))=c_F
\end{equation}
for some constant $c_F>0$ and some $3<\tau <4$, then there is a constant $b>0$ such that for all $A>1$ and $n\geq 1$, the $NR_n(\mathbf{w})$ random graph satisfies
\begin{equation}\label{pp2}
\mathbb{P}(|\mathcal{C}_{\max}|\notin [A^{-1}n^{(\tau-2)/(\tau-1)},An^{(\tau-2)/(\tau-1)}])\leq b/A.
\end{equation}
On the other hand, when
\begin{equation}\label{six}
1-F(x)\leq c_Fx^{-(\tau-1)}\hspace{0.2cm}(x\geq 0)
\end{equation}
for some $c_F>0$ and some $\tau>4$, then there is a constant $b>0$ such that, for all $n\geq 1$ and all $A>1$, the $NR_n(\mathbf{w})$ random graph satisfies
\begin{equation}\label{pp1}
\mathbb{P}(|\mathcal{C}_{\max}|\notin [A^{-1}n^{2/3},An^{2/3}])\leq b/A.
\end{equation}
(Actually Van der Hofstad \cite{hofstad_critical_behaviour} proved a more general result, namely that the lower bounds (\ref{pp2}) and (\ref{pp1}) remain valid also after a \textit{small perturbation} of the vertex weights; see Theorems 1.1 and 1.2 in \cite{hofstad_critical_behaviour}.) 

For an heuristic explanation of the critical behaviour described by (\ref{pp2}) and (\ref{pp1}), we refer to section 1.3 in \cite{hofstad_critical_behaviour}.

We also mention that in \cite{de_ambroggio_pachon:upper_bounds_inhom_RGs} the authors used the first part of the martingale argument introduced by Nachmias and Peres \cite{nachmias_peres:CRG_mgs} to obtain simple upper bounds for the probability of observing unusually large maximal components in the (critical) $NR_n(\mathbf{w})$ random graph for both regimes $\tau\in (3,4)$ and $\tau>4$, even if in the former case (i.e. for $\tau\in (3,4)$) the distribution function $F$ is required to satisfy a stronger condition with respect to the one stated in (\ref{five}).

Our goal here is to use Theorem \ref{Umbi} to provide a very simple proof of the fact that, in the critical $NR_n(\mathbf{w})$ model with vertex weights as in (\ref{wj}) and distribution function $F$ satisfying (\ref{six}), the largest component is unlikely to contain more than $n^{2/3}$ vertices. More precisely, we prove the following
\begin{prop}\label{CIRG}
	Consider the $NR_n(\mathbf{w})$ random graph with weights defined as in (\ref{wj}) above. Suppose that there exists a constant $c_F>0$ and a $\tau >4$ such that $1-F(x)\leq c_Fx^{-(\tau-1)}$ for all $x\geq 0$. Then, given any $A>1$, when $n$ is large enough we have that
	\begin{align*}
	\mathbb{P}(|\mathcal{C}_{\max}|> \lceil An^{2/3} \rceil)\leq \frac{c_3}{A^{3/2}},
	\end{align*}
	where $c_3$ is a finite positive constant which depends solely on $c_F$ and $\tau$.
\end{prop}

\section{Proofs}
Here we are going to prove the results stated in Section 2. We start by proving Theorem \ref{Umbi} and subsequently we prove the remaining results, namely Propositions \ref{CRIG}, \ref{CBPRRG} and \ref{CIRG}.

The proof of Theorem \ref{Umbi} relies on the following \textit{ballot-type} estimate, which is taken from \cite{de_ambroggio_roberts1}. For a general introduction to classical ballot theorems and their generalisations, see for instance \cite{addario_berry_reed:ballot_theorems}, \cite{kager:hitting_time} and references therein. 
\begin{lem}\label{ballot}
	Fix $n\in \mathbb{N}$ and let $(W_i)_{i\geq 1}$ be a sequence of i.i.d.~valued random variables taking values in $\mathbb{Z}$. Let $r\in \mathbb{N}$, and suppose that $\mathbb{P}(W_1=r)>0$. Define $S_t = \sum_{i=1}^{t}W_i$ for $t\in \mathbb{N}_0$. Then, for any $j\geq 1$, we have 
	\begin{equation*}
	\mathbb{P}(r+S_t>0\hspace{0.15cm} \forall t\in [n],r+S_{n}=j)\leq \mathbb{P}(X_1=r)^{-1}\frac{j}{n+1}\mathbb{P}(S_{n+1}=j).
	\end{equation*}
\end{lem}
\subsection{Proof of Theorem \ref{Umbi}.}
Let $k=k(n)\in \mathbb{N}$. By hypothesis, there is a sequence of i.i.d.~random variables $X_i$ taking values in $\mathbb{N}_0$ such that, setting $S_t= \sum_{i=1}^{t}(X_i-1)$, then
\begin{equation*}
\mathbb{P}(|\mathcal{C}(V_n)|>k)\leq \mathbb{P}(2+S_t>0 \hspace{0.15cm}\forall t\in [k]).
\end{equation*}
Using Lemma \ref{ballot} with $W_i=X_i-1$ and $r=2$ we obtain
\begin{align}\label{prima}
\nonumber\mathbb{P}(2+S_t>0 \hspace{0.15cm}\forall t\in [k])&= \sum_{h=1}^{\infty}\mathbb{P}(2+S_t>0 \hspace{0.15cm}\forall t\in [k], 2+S_k=h)\\
&\leq a\sum_{h=1}^{\infty}\frac{h}{k+1}\mathbb{P}(S_{k+1}=h),
\end{align}
where we set $a= 1/\mathbb{P}(X_1=3)\in (0,\infty)$. Now let $m$ be a non-negative integer to be specified later. By splitting the series in (\ref{prima}) at $h=m$ we can write
\begin{equation}\label{sec}
a\sum_{h=1}^{\infty}\frac{h}{k+1}\mathbb{P}(S_{k+1}=h)\leq a\frac{m}{k+1}+\frac{a}{k+1}\sum_{h=m+1}^{\infty}h\mathbb{P}(S_{k+1}=h).
\end{equation}
Now the series in (\ref{sec}) equals
\begin{align}\label{four}
\nonumber\frac{a}{k+1}&\sum_{h=m+1}^{\infty}h\mathbb{P}\left(\sum_{i=1}^{k+1}X_i=h+k+1\right)\\
&=\frac{a}{k+1}\sum_{z=m+k+2}^{\infty}z\mathbb{P}\left(\sum_{i=1}^{k+1}X_i=z\right)-a\mathbb{P}\left(\sum_{i=1}^{k+1}X_i\geq m+k+2\right).
\end{align}
To proceed, we observe the following: if $X$ is a random variable taking values in $\mathbb{N}_0$, then for any $h\ge1$, we have
\begin{align*}
\mathbb{E}\left(X\mathbbm{1}_{\{X\geq h\}}\right) = \mathbb{E}\left(\sum_{i=1}^{\infty} \mathbbm{1}_{\{i\le X\}}\mathbbm{1}_{\{X\geq h\}}\right)&= h\mathbb{P}(X\geq h) + \sum_{i=h+1}^{\infty} \mathbb{P}(X\geq i).
\end{align*}
Thus the series in (\ref{four}) equals
\begin{align*}
a\frac{m+1}{k+1}\mathbb{P}\left(\sum_{i=1}^{k+1}X_i\geq m+k+2\right)&+a\mathbb{P}\left(\sum_{i=1}^{k+1}X_i\geq m+k+2\right)\\
&+\frac{a}{k+1}\sum_{z=m+k+3}^{\infty}\mathbb{P}\left(\sum_{i=1}^{k+1}X_i\geq z\right).
\end{align*}
Substituting the series in (\ref{four}) with these three terms we obtain
\begin{multline}\label{sumsum}
\frac{a}{k+1}\sum_{h=m+1}^{\infty}h\mathbb{P}(S_{k+1}=h)\\
= a\frac{m}{k+1}\mathbb{P}\left(\sum_{i=1}^{k+1}X_i\geq m+k+2\right)+\frac{a}{k+1}\sum_{z=m+k+2}^{\infty}\mathbb{P}\left(\sum_{i=1}^{k+1}X_i\geq z\right).
\end{multline}
Now observe that the series in (\ref{sumsum}) can be rewritten as follows:
\begin{equation*}
\frac{a}{k+1}\sum_{z=m+k+2}^{\infty}\mathbb{P}\left(\sum_{i=1}^{k+1}X_i\geq z\right)=\frac{a}{k+1}\sum_{h=m+1}^{\infty}\mathbb{P}\left(\sum_{i=1}^{k+1}X_i\geq h+k+1\right).
\end{equation*}
Summarizing, so far we have shown that
\begin{multline}\label{cinque}
\mathbb{P}(|\mathcal{C}(V_n)|>k)\leq a\frac{m}{k+1}+a\frac{m}{k+1}\mathbb{P}\left(\sum_{i=1}^{k+1}X_i\geq k+1+(m+1)\right)\\
+\frac{a}{k+1}\sum_{h=m+1}^{\infty}\mathbb{P}\left(\sum_{i=1}^{k+1}X_i\geq k+1+h\right).
\end{multline}
Using or assumption on $\mathbb{E}(e^{rX_1})$ and Markov's inequality we have that, for all $h\geq m+1$ and all $r\in (0,\rho)$, 
\begin{align*}
\mathbb{P}\left(\sum_{i=1}^{k+1}X_i\geq k+1+h\right)&\leq e^{-r(k+1)-rh}\mathbb{E}\left(e^{rX_1}\right)^{k+1}\\
&\leq \exp\left\{-r(k+1)-rh+r(1+\epsilon)(k+1)+\delta r^2 (k+1)\right\}\\
&=\exp\left\{-rh+r\epsilon (k+1)+\delta r^2 (k+1)\right\}.
\end{align*}
Now if $k$ is such that $\rho \sqrt{k}\geq 1$ for all large enough $n$, then $r= 1/\sqrt{k+1}<\rho$ and hence using this specific value of $r$ we obtain
\begin{equation*}
\mathbb{P}\left(\sum_{i=1}^{k+1}X_i\geq k+1+h\right)\leq e^{-\frac{h}{\sqrt{k+1}}}e^{\epsilon \sqrt{k+1}+\delta}.
\end{equation*}
Also, if $k$ satisfies $\epsilon \sqrt{k}\leq c$ for all sufficiently large $n$, we see that $\epsilon \sqrt{k+1}\leq 2c$ and hence we can bound from above the series in (\ref{cinque}) as follows:
\begin{equation}\label{eleven}
\frac{a}{k+1}\sum_{h=m+1}^{\infty}\mathbb{P}\left(\sum_{i=1}^{k+1}X_i\geq h+k+1\right)\leq \frac{ae^{\delta+2c}}{k+1}\sum_{h=m+1}^{\infty}e^{-\frac{h}{\sqrt{k+1}}}.
\end{equation}
Now observe that
\begin{equation*}
\sum_{h=m+1}^{\infty}e^{-\frac{h}{\sqrt{k+1}}}=e^{-\frac{m+1}{\sqrt{k+1}}}\frac{1}{1-e^{-\frac{1}{\sqrt{k+1}}}}.
\end{equation*}
Using the inequality $e^{-x}\leq 1-x+x^2/2$ (which is valid for all $x\geq 0$) we see that
\begin{equation*}
1-e^{-\frac{1}{\sqrt{k+1}}}\geq \frac{1}{\sqrt{k+1}}-\frac{1}{2(k+1)}\geq \frac{1}{2\sqrt{k+1}},
\end{equation*}
and hence the expression on the right-hand side of (\ref{eleven}) is at most
\begin{equation*}
\frac{ae^{\delta+2c}}{k+1}e^{-\frac{m+1}{\sqrt{k+1}}}2\sqrt{k+1}=\frac{2ae^{\delta+2c}}{\sqrt{k+1}}e^{-\frac{m+1}{\sqrt{k+1}}}.
\end{equation*}
Thus we obtain
\begin{equation*}
\mathbb{P}(|\mathcal{C}(V_n)|>k)\leq a\frac{m}{k+1}+ae^{\delta+2c}\frac{m}{k+1} e^{-\frac{m+1}{\sqrt{k+1}}}+\frac{2ae^{\delta+2c}}{\sqrt{k+1}}e^{-\frac{m+1}{\sqrt{k+1}}}.
\end{equation*}
Taking $m=\lfloor \sqrt{k+1} \rfloor$ we see that
\begin{equation*}
\mathbb{P}\left(|\mathcal{C}(V_n)|>k\right)\leq \frac{a}{\sqrt{k+1}}+3ae^{\delta+2c}\frac{e^{-1}}{\sqrt{k+1}}=C\frac{a}{\sqrt{k+1}},
\end{equation*}
where we set $C=C(\delta,c)= 1+3e^{\delta+2c-1}$. Finally, denoting by
\begin{equation*}
N_k=\sum_{v\in [n]}^{}\mathbbm{1}_{\{|\mathcal{C}(v)|> k \}}
\end{equation*} 
the number of vertices that are contained in components of size at least $k$ we obtain
\begin{equation*}
\mathbb{P}(|\mathcal{C}_{\max}|> k)= \mathbb{P}(N_{k}> k)\leq \frac{1}{k}\mathbb{E}(N_{k})=\frac{n}{k}\mathbb{P}(|\mathcal{C}(V_n)|> k)\leq \frac{C}{\mathbb{P}(X_1=3)}\frac{n}{k^{3/2}},
\end{equation*}
completing the proof of the theorem.

\subsection{Proof of Proposition \ref{CRIG}.}
Let $\mathbb{H}(n,m,p)$ be a random (multi-)graph constructed from the bipartite graph $\mathbb{B}(n,m,p)$ by letting the number of edges between $v_i,v_j\in V$ equal the number of auxiliary vertices $a_k$ that are adjacent to both $v_i$ and $v_j$. (Recall that $V$ is the vertex set of the random intersection graph $\mathbb{G}(n,m,p)$ under investigation.) Notice that $\mathbb{G}(n,m,p)$ can be obtained from $\mathbb{H}(n,m,p)$ by coalescing multiple edges between vertices into one single edge. Hence, thanks to this construction, we see that the degree distribution in $\mathbb{G}(n,m,p)$ is dominated by the degree distribution in $\mathbb{H}(n,m,p)$. Notice that the latter is a compound binomial distribution with moment generating function
\begin{equation}\label{mgf}
\mathbb{E}\left(e^{rX_1}\right)=\left\{1-\frac{\gamma}{n}+\frac{\gamma}{n}\left(1-\frac{\gamma}{n}+\frac{\gamma}{n}e^{r}\right)^{n-1}\right\}^{\lfloor \beta n \rfloor},
\end{equation} 
since (by construction) a vertex $v\in \mathbb{H}(n,m,p)$ is connected to a $\text{Bin}(m,p)$ number of auxiliary vertices, and each one of them is connected to an independent $\text{Bin}(n-1,p)$ number of vertices in $V\setminus \{v\}$.

Therefore, by revealing the components of $\mathbb{G}(n,m,p)$ using the exploration process described at the beginning of Section 2, we can write 
\begin{align*}
\mathbb{P}(|\mathcal{C}(V_n)|>k)\leq \mathbb{P}\left(2+\sum_{i=1}^{t}(X_i-1)>0\hspace{0.15cm}\forall t\in [k]\right),
\end{align*}
where the $X_i$ are i.i.d.~compound binomial random variables with moment generating function given in (\ref{mgf}).

Using the probability generating function of $X_1$ (which coincides with (\ref{mgf}) after substituting $e^r$ with $r$), it is not difficult to show that (for large enough $n$) the probability $\mathbb{P}(X_1=3)$ is bounded from below by a positive constant which depends solely on $\gamma$ and $\beta$. 

Next, in order to apply Theorem \ref{Umbi}, we simply need to prove an upper bound for $\mathbb{E}\left(e^{rX_1}\right)$. Recalling the expression of the moment generating function of $X_1$ given in (\ref{mgf}) we obtain
\begin{equation}\label{long}
\mathbb{E}\left(e^{rX_1}\right)=\exp\left\{\lfloor \beta n \rfloor \log\left(1-\frac{\gamma}{n}+\frac{\gamma}{n}\left[1+\frac{\gamma}{n}\left(e^{r}-1\right)\right]^{n-1}\right)\right\}.
\end{equation}
Taking $r\in (0,1)$ we have that $e^r-1\leq r+r^2$. Then, since $1+x\leq e^x$ for all $x\in \mathbb{R}$, we see that (\ref{long}) is at most
\begin{equation*}\label{hhhh}
\exp\left\{\lfloor \beta n \rfloor \log\left(1+\frac{\gamma}{n}\left(\exp\left\{\gamma(r+r^2)\right\}-1\right)\right)\right\}.
\end{equation*}
Thus taking $r<\min\{1,1/2\gamma\}$ (so that $\gamma(r+r^2)<1$) we see that
\begin{equation*}
\exp\left\{\gamma(r+r^2)\right\} - 1\leq \gamma(r+r^2)+\gamma^2(r+r^2)^2
\end{equation*}
and hence, using the fact that $\log(1+x)\leq x$ for all $x>-1$, we can write
\begin{multline}\label{fffff}
\exp\left\{\lfloor \beta n \rfloor \log\left(1+\frac{\gamma}{n}\left(\exp\left\{\gamma(r+r^2)\right\}-1\right)\right)\right\}\\ 
\leq \exp\left\{\lfloor \beta n \rfloor\left(\frac{\gamma^2}{n}(r+r^2) + \frac{\gamma^3}{n}(r+r^2)^2\right)\right\}.
\end{multline}
Recalling that $\beta \gamma^2=1$ we obtain
\begin{equation*}
\exp\left\{\lfloor \beta n \rfloor\left(\frac{\gamma^2}{n}(r+r^2) + \frac{\gamma^3}{n}(r+r^2)^2\right)\right\}\leq  \exp \left\{r+r^2(1+4\gamma)\right\}.
\end{equation*}
Therefore, for all $r\in (0,\min\{1,1/2\gamma\})$ we have that
\begin{equation*}
\mathbb{E}(e^{rX_1})\leq e^{r+r^2(1+4\gamma)}
\end{equation*}
and hence condition $(iii)$ in Theorem \ref{Umbi} is satisfied for $\rho=\min\{1,1/2\gamma\}$, $\delta=1+4\gamma$ and $\epsilon=0$.
Hence, taking $k=\lceil A n^{2/3} \rceil$ (which clearly satisfies the requirement $\epsilon \sqrt{k}\leq c$ for $c=0$ as well as the condition $\rho \sqrt{k}\geq 1$) and applying Theorem \ref{Umbi} we obtain that
\begin{align*}
\mathbb{P}(|\mathcal{C}_{\max}|>\lceil A n^{2/3} \rceil)\leq \frac{c_1}{A^{3/2}},
\end{align*}
for some finite positive constant $c_1$ which depends solely on $\gamma$ and $\beta$.

\subsection{Proof of Proposition \ref{CBPRRG}.}
Since $\mathbb{G}$ is $d$-regular, we can use the exploration process described at the beginning of Section 2 to conclude that 
\begin{align*}
\mathbb{P}(|\mathcal{C}(V_n)|>k)\leq \mathbb{P}\left(2+\sum_{i=1}^{t}(X_i-1)>0\hspace{0.15cm}\forall t\in [k]\right),
\end{align*}
where the $X_i$ are i.i.d.~random variables with $X_i\sim Bin(d-1,p)$, so that $\sum_{i=1}^{k+1}X_i\sim \text{Bin}((k+1)(d-1),p)$. (Notice that, if we would have started the exploration process with only one active vertex, now we would have $\eta_1\sim \text{Bin}(d,p)$ and hence in particular it would be impossible to dominate $\eta_1$ with a $\text{Bin}(d-1,p)$ random variable.) Using a monotonicity argument we can focus on the (critical) case $p=1/(d-1)$. Note that, since $d>3$,
\begin{align*}
\mathbb{P}(X_1=3)=\frac{(d-2)(d-3)}{6(d-1)^2}\left(1-\frac{1}{d-1}\right)^{d-4}>0.
\end{align*}
Next, for all $r\in (0,1)$, using the inequality $1+x\leq e^x$ (valid for all $x\in \mathbb{R}$) we see that
\begin{align}\label{binbound}
\mathbb{E}\left(e^{rX_1}\right)=\left(1+\frac{1}{d-1}(e^r-1)\right)^{d-1}\leq e^{e^r-1}\leq e^{r+r^2}
\end{align}
Hence condition $(iii)$ of Theorem \ref{Umbi} is satisfied for $\rho,\delta=1$ and $\epsilon=0$. Thus taking $k=\lceil An^{2/3}\rceil$ (which satisfies $\epsilon \sqrt{k}\leq c$ for $c=0$, as well as $\rho \sqrt{k}\geq 1$) we arrive at
\begin{align*}
\mathbb{P}(|\mathcal{C}_{\max}|>k)\leq c_2\frac{n}{k^{3/2}},
\end{align*}
for some finite positive constant $c_2$ which depends solely on $d$.

\subsection{Proof of Proposition \ref{CIRG}.}
Before starting with the actual proof, we need to recall the definition of \textit{size-biased} distribution of a non-negative random variable and to introduce a few facts. 
\begin{defn}
	For a non-negative random variable $X$ with $\mathbb{E}(X)>0$, the size-biased distribution of $X$, denote by $X^*$, is the random variable defined by
	\begin{equation*}
	\mathbb{P}(X^*\leq x)=\frac{\mathbb{E}\left(X\mathbbm{1}_{\{X\leq x\}}\right)}{\mathbb{E}(X)}.
	\end{equation*} 
\end{defn}
For proofs of the assertions that appear in the statement of next result, see \cite{de_ambroggio_pachon:upper_bounds_inhom_RGs}.
\begin{lem}\label{summmmmm}
	Suppose that $1-F(x)\leq c_Fx^{-(\tau-1)}$ for all $x\geq 0$, for some $c_F>0$ and $\tau>4$. Let $w_i$ be as in (\ref{wj}). Then $\max\{w_i:i\in [n]\}\leq (c_F n)^{1/(\tau-1)}$. Moreover, defining 
	\begin{equation}\label{fn}
	F_n(x)= \frac{1}{n}\sum_{i=1}^{n}\mathbbm{1}_{\{w_i\leq x\}}
	\end{equation}
	and letting $W_n$ being a random variable with distribution function $F_n$ and size-biased distribution $W_n^*$, then $\mathbb{E}\left((W^*_n)^2\right)\leq C_1$ and $|1-\mathbb{E}\left(W^*_n\right)|\leq C_2n^{-\frac{\tau-3}{\tau-1}}$ for all large enough $n$, where $C_1$ and $C_2$ are two positive constants which depend on $c_F$ and $\tau$.
\end{lem}
As explained in Van der Hofstad \cite{hofstad_critical_behaviour} (see also \cite{de_ambroggio_pachon:upper_bounds_inhom_RGs}), the cluster exploration of $V_n$ in the $NR_n(\mathbf{w})$ random graph can be dominated by the total progeny of a (marked mixed-Poisson) branching process. Specifically, following Van der Hofstad \cite{hofstad_critical_behaviour} we can write
\begin{equation}
\mathbb{P}\left(|\mathcal{C}(V_n)|>k\right)\leq \mathbb{P}\left(2+\sum_{i=1}^{t}(X_i-1)>0\hspace{0.15cm}\forall t\in [k]\right),
\end{equation}
where the $X_i$ are independent mixed Poisson random variables with $X_i\sim \text{Poi}(w_{M_i})$ and the $M_i$ are i.i.d.~random variables, all distributed as a random variable $M$ with distribution given by
\[\mathbb{P}(M=m)=\frac{w_m}{l_n}, \hspace{0.2cm}m\in [n].\]
As remarked in \cite{hofstad_critical_behaviour}, a $\text{Poi}(w_M)$ random variable converges in distribution to a mixed Poisson random variable with random parameter $W^*$, where $W^*$ is the size-biased distribution of $W$, the latter being a positive random variable with distribution function $F$. Therefore, $\mathbb{P}(X_1=3)$ converges to $\mathbb{P}(Z=3)$, where $Z\sim \text{Poi}(W^*)$. It follows that $\mathbb{P}(X_1=3)\geq \mathbb{P}(Z=3)/2$ for all large enough $n$, and hence we obtain
\begin{equation}\label{p3}
\mathbb{P}(X_1=3)\geq\frac{1}{2}\mathbb{E}\left(\mathbb{P}(Z=3|W^*)\right)=\mathbb{E}\left(e^{-W^*}(W^*)^3\right)/12>0.
\end{equation}
Taking $r\in (0,1)$ so that $e^r-1\leq r+r^2$ we obtain
\begin{multline*}
\mathbb{E}\left[e^{rX_1}\right]=\sum_{h=0}^{\infty}e^{rh}\sum_{i\in [n]}^{}\frac{w_i}{l_n}e^{-w_i}\frac{w^{h}_i}{h!}=\sum_{i\in [n]}^{}\frac{w_i}{l_n}e^{w_i(e^r-1)}\\
\leq \sum_{i\in [n]}^{}\frac{w_i}{l_n}e^{w_i(r+r^2)}=\exp\left\{\log\left(\sum_{i\in [n]}^{}\frac{w_i}{l_n}e^{w_i(e^r-1)}\right)\right\}.
\end{multline*}
If moreover $r<(2\max\{w_i:i\in [n]\})^{-1}$ then $w_i(r+r^2)<2r w_i\leq 1$ for every $i\in [n]$ and hence, for $0<r<\min\{1,1/2\max\{w_i:i\in [n]\}\}$, we can bound
\begin{multline*}
\log\left(\sum_{i\in [n]}^{}\frac{w_i}{l_n}e^{w_i(e^r-1)}\right)\leq \log \left(\sum_{i\in [n]}^{}\frac{w_i}{l_n}(1+w_i(r+r^2)+w^2_i4r^2)\right)\\
=\log\left(1+(r+r^2)\sum_{i\in [n]}^{}\frac{w^2_i}{l_n}+4r^2\sum_{i\in [n]}^{}\frac{w^3_i}{l_n}\right).
\end{multline*}
Note that, if $W_n$ is a random variable with distribution function $F_n$ given in (\ref{fn}) and $W_n^*$ is its size-biased distribution, then
\begin{equation*}
\sum_{i\in [n]}^{}\frac{w^2_i}{l_n}=\mathbb{E}(W_n^*) \text{\hspace{0.2cm} and \hspace{0.2cm}}\sum_{i\in [n]}^{}\frac{w^3_i}{l_n}=\mathbb{E}\left((W^*_n)^2\right).
\end{equation*}
Therefore we arrive at
\begin{equation*}
\log\left(\sum_{i\in [n]}^{}\frac{w_i}{l_n}e^{w_i(e^r-1)}\right)\leq \log\left(1+(r+r^2)\mathbb{E}(W_n^*)+4r^2\mathbb{E}\left((W^*_n)^2\right)\right).
\end{equation*}
Thanks to Lemma \ref{summmmmm} we know that $\mathbb{E}\left((W^*_n)^2\right)\leq C_1$ and $\mathbb{E}\left(W^*_n\right)\leq 1+C_2n^{-\frac{\tau-3}{\tau-1}}$ for all sufficiently large $n$ (for some finite constants $C_1,C_2>0$ which depend on $c_F$ and $\tau$), we see that
\begin{equation*}
\log\left(1+(r+r^2)\mathbb{E}(W_n^*)+4r^2\mathbb{E}\left((W^*_n)^2\right)\right)\leq r\left(1+C_2n^{-\frac{\tau-3}{\tau-1}}\right)+r^2(1+5C_1)
\end{equation*}
for all large enough $n$. Summarizing, for all positive $r<\min\{1,1/2\max\{w_i:i\in [n]\}\}$ we have shown that
\begin{align*}
\mathbb{E}\left(e^{rX_1}\right)\leq \exp\left\{r\left(1+C_2n^{-\frac{\tau-3}{\tau-1}}\right)+r^2(1+5C_1)\right\},
\end{align*}
provided $n$ is sufficiently large. Hence condition $(iii)$ in Theorem \ref{Umbi} is satisfied for $\rho=\rho(n)=\min\{1,1/2\max\{w_i:i\in [n]\}\}$, $\delta=1+5C_1$ and  $\epsilon=\epsilon(n)=C_2n^{-\frac{\tau-3}{\tau-1}}$. Note that $k= \lceil An^{2/3} \rceil$ satisfies $\epsilon \sqrt{k}=O(A^{1/2}C_2n^{-\frac{\tau-4}{3(\tau-1)}})\leq 1$ for all large enough $n$ since $A$ is fixed and $\tau>4$, and $\rho \sqrt{k}\geq 1$ since $(2\max\{w_i:i\in [n]\})^{-1}\sqrt{k}\geq A^{1/2}n^{\frac{\tau-4}{3(\tau-1)}}(2c^{1/(\tau-1)}_F)^{-1}\geq 1$. Whence we can apply Theorem \ref{Umbi} to conclude that
\begin{align*}
\mathbb{P}(|\mathcal{C}_{\max}|>\lceil An^{2/3} \rceil)\leq \frac{c_3}{A^{3/2}},
\end{align*}
for some finite positive constant $c_3$ that depends solely on $c_F$ and $\tau$.\\

\textbf{Acknowledgements.} The author would like to thank the Royal Society for his PhD scholarship and M. Roberts for useful suggestions that helped improving the presentation of the paper. Moreover, the author thanks G. Perarnau and A. Pachon for interesting discussions about random intersection graphs in occasion of the workshop \textit{Graphs and Randomness in Turin} (January 2019). \\

\bibliographystyle{plain}

\begin{thebibliography}{99}
	\bibliographystyle{plain}
	
	\bibitem{addario_berry_reed:ballot_theorems}
	{\sc Addario-Berry, L., Bruce, R.} (2008). Ballot theorems, old and new. {\em Horizons of Combinatorics}
	
	\bibitem{beh2007}
	{\sc Behrisch, M.} (2007). Component evolution in random intersection graphs. {\em The Electronic Journal of Combinatorics}.
	
	
	\bibitem{bollobas_book}
	{\sc Bollob{\'a}s, B.} (2001). {\em Random graphs}, 2nd~edn. Cambridge University Press, Cambridge.
	
	\bibitem{boll_janson_riord_inhom}
	{\sc Bollobas, B., Janson, S., Riordan, O.} (2007). The phase transition in inhomogeneous random graphs. {\em Random Structures \& Algorithms}.
	
	
	\bibitem{chung2002connected}
	{\sc Chung, F., Lu, L.} (2002). Connected components in random graphs with given expected degree sequences. {\em Annals of Combinatorics}.
	
	\bibitem{chung2004average}
	{\sc Chung, F., Lu, L.} (2003). The average distance in a random graph with given expected degrees.  {\em Internet Mathematics}.
	
	
	\bibitem{chung2006volume}
	{\sc Chung, F., Lu, L.} (2006). The volume of the giant component of a rand om graph with given expected degrees. {\em SIAM Journal of Discrete Mathematics}.
	
	
	\bibitem{de_ambroggio_roberts1}
	{\sc De Ambroggio, U., Roberts, M.~I.} (2021). Unusually large components in near-critical Erd\H{o}s-R\'enyi graphs via ballot theorems. \newblock Preprint:  \texttt{http://arxiv.org/abs/2101.05358}.
	
	\bibitem{de_ambroggio_pachon:upper_bounds_inhom_RGs}
	{\sc De Ambroggio, U., Pachon, A.} (2020). Simple upper bounds for the largest components in critical inhomogeneous random graphs. \newblock Preprint:  \texttt{http://arxiv.org/abs/2012.09001}.
	
	\bibitem{deijfen_kets_2009}
	{\sc Deijfen, M., Kets, W.} (2009). Random intersection graphs with tunable degree distribution and clustering. {\em Probability in the Engineering and Informational Sciences}.
	
	\bibitem{frieze_karo?ski_2015}
	{\sc Frieze, A., Karonski, M.} (2015). {\em Introduction to random graphs}, Cambridge University Press, Cambridge.
	
	\bibitem{remco:random_graphs}
	{\sc van der Hofstad, R.} (2016). {\em Random graphs and complex networks (Volume 1)}, Cambridge University Press, Cambridge.
	
	\bibitem{hofstad_critical_behaviour}
	{\sc van der Hofstad, R.} (2013). Critical behaviour in inhomogeneous random graphs. {\em  Random Structures \& Algorithms}.
	
	\bibitem{janson_et_al:random_graphs}
	{\sc Janson, S., Luczak, T., Rucinski, A.} (2011). {\em Random graphs}, John Wiley \& Sons, New York.
	
	\bibitem{janson_contiguity}
	{\sc Janson, S.} (2009). Asymptotic equivalence and contiguity of some random graphs. {\em Random Structures \& Algorithms}.
	
	\bibitem{joosper}
	{\sc Joos, F., Perarnau, G.} (2018). Critical percolation on random regular graphs. {\em Proceedings of the American Mathematical Society}.
	
	\bibitem{kager:hitting_time}
	{\sc Kager, W.} (2011). The hitting time theorem revisited. {\em The American Mathematical Monthly}.
	
	\bibitem{kang2018evolution}
	{\sc Kang, M., Pachon, A., Rodriguez, P.M.} (2018). Evolution of a Modified Binomial Random Graph by Agglomeration. {\em Journal of Statistical Physics}.
	
	\bibitem{KrivSud}
	{\sc Krivelevich, M., Sudakov, B.} (2012). The phase transition in random graphs: A simple proof. {\em Random Structures \& Algorithms}.
	
	
	\bibitem{laglind2008}
	{\sc Lageras, A.N., Lindholm, M.} (2008). A note on the component structure in random intersection graphs with tunable clustering. {\em The Electronic Journal of Combinatorics}.
	
	
	\bibitem{luczak_et_al:structure_RG}
	{\sc Luczak, T., Pittel, B., Wierman, J.~C.} (1994). The structure of a random graph at the point of the phase transition. {\em Transactions of the American Mathematical Society}.
	
	
	
	
	\bibitem{nachmiasperes_criticalperc}
	{\sc Nachmias,A., Peres, Y.} (2010). Critical percolation on random regular graphs. {\em Random Structures \& Algorithms}.
	
	\bibitem{nachmias_peres:CRG_mgs}
	{\sc Nachmias,A., Peres, Y.} (2010). The critical random graph, with martingales. {\em Israel Journal of Mathematics}.
	
	
	
	
	\bibitem{norros_reittu_2006}
	{\sc Norros, I., Reittu, H.} (2006). On a conditionally Poissonian graph process. {\em Advances in Applied Probability}.
	
	
	
	
	\bibitem{penrose_2018}
	{\sc Penrose, M.D.} (2018). Inhomogeneous random graphs, isolated vertices, and Poisson approximation. {\em Journal of Applied Probability}.
	
	\bibitem{pittel:largest_cpt_rg}
	{\sc Pittel, B.} (2001). On the largest component of the random graph at a nearcritical stage. {\em Journal of Combinatorial Theory, Series B}.
	
	
	\bibitem{roberts:component_ER}
	{\sc Roberts, M.~I.} (2017). The probability of unusually large components in the near-critical {E}rd{\H{o}}s-{R}\'{e}nyi graph. {\em Advances in Applied Probability}.
	
	\bibitem{starkRIG}
	{\sc Stark, D.} (2004). The vertex degree distribution of random intersection graphs. {\em Random Structures and Algorithms}.
	
	
	
\end{thebibliography}

\end{document}